\documentclass[12pt,reqno]{amsart}
\usepackage{amsmath,amssymb,amsfonts,amsthm,amstext}
\usepackage{url,xspace}
\usepackage{color,graphicx}
\usepackage{graphics}
\usepackage{enumerate}
\usepackage{cleveref}
\usepackage{tikz}
\usepackage[margin=3cm]{geometry}

\renewcommand{\P}{\mathbb{P}}

\newcommand{\ind}{{\mathbf 1}}
\newcommand{\id}{\mbox{\rm id}}

\renewcommand{\SS}{\mathcal{S}}

\newcommand{\df}{\textbf}

\DeclareMathOperator{\rank}{rank}

\newcommand{\dia}{{\hfill $\Diamond$}}

\newtheorem{thm}{Theorem}
\newtheorem{lemma}[thm]{Lemma}
\newtheorem{prop}[thm]{Proposition}

\newtheorem{question}{Question}

\theoremstyle{definition}
\newtheorem{construction}{Construction}

\crefname{thm}{Theorem}{Theorems}
\crefname{lemma}{Lemma}{Lemmas}
\crefname{prop}{Proposition}{Propositions}
\crefname{cor}{Corollary}{Corollaries}
\crefname{section}{Section}{Sections}
\crefname{figure}{Figure}{Figures}
\crefname{table}{Table}{Tables}
\crefname{definition}{Definition}{Definitions}
\crefname{construction}{Construction}{Constructions}

\newcounter{icount}
\newenvironment{ilist}{\begin{list}{{\rm (\roman{icount})}}%
{\usecounter{icount}\itemsep 1pt}}{\end{list}}

\title{Perfect shuffling by lazy swaps}
\author{Omer Angel}
\address{Omer Angel, Mathematics, University of British Columbia, Canada}
\email{angel@math.ubc.ca}
\author{Alexander E. Holroyd}
\address{Alexander E. Holroyd}
\email{holroyd@uw.edu}

\date{20 February 2018}

\keywords{Random permutation; shuffling; transposition; reduced word; perfect mixing; sorting network}
\subjclass[2010]{05A05; 60C05; 68P10}

\begin{document}
\maketitle

\begin{abstract}
We characterize the minimum-length sequences of independent lazy simple transpositions whose composition is a uniformly random permutation.  For every reduced word of the reverse permutation there is exactly one valid way to assign probabilities to the transpositions.  It is an open problem to determine the minimum length of such a sequence when the simplicity condition is dropped.
\end{abstract}

\section{Introduction}

Let $S_n$ be the symmetric group of all permutations of $1,\ldots,n$, with composition given by $(\sigma\tau)(i):=\sigma(\tau(i))$ for $1\leq i\leq n$.  A \df{lazy transposition} with parameters $(a,b,p)$ is a random permutation $T$ that with probability $p$ equals the transposition (or swap) $t(a,b):=(a\; b)\in S_n$ exchanging the elements in positions $a$ and $b$, and otherwise equals the identity $\id\in S_n$.  Given a sequence of parameters $\SS=(a_i,b_i,p_i)_{i=1}^\ell$, let $T_1,\ldots,T_\ell$ be independent lazy transpositions, where $T_i$ has parameters $(a_i,b_i,p_i)$.  We say that $\SS$ is a (perfect) \df{transposition shuffle} (of order $n$ and length $\ell$) if the composition $T_1\cdots T_\ell$ of these random permutations is uniformly distributed on $S_n$.  We pose the following apparently unsolved question.

\begin{question}
What is the minimum possible length $L_n$ of a transposition shuffle of order $n$?  Is it the case that $L_n={n\choose 2}$ for all $n$?
\end{question}

The best bounds we know for general $n$ are
$$ \log_2 n! \leq L_n \leq {n\choose 2}.$$
The lower bound (which is of course asymptotic to $n\log_2 n$ as $n\to\infty$) follows by the obvious counting argument: a composition of $\ell$ lazy transpositions can take at most $2^\ell$ possible values, while $\#S_n=n!$.  In the other direction, we have several distinct constructions of transposition shuffles of length exactly $n\choose 2$ (for all $n$), and none shorter (for any $n$).  It can be verified by case analysis that $L_n={n \choose 2}$ for $n\leq 4$.  Computer experiments by Viktor Kiss (personal communication) suggest that the same holds for $n=5$ also.

Our main result addresses the special case
of \df{simple} transposition shuffles, by which we mean those that transpose only adjacent pairs: $b_i=a_i+1$ for all $i$.  We will characterize the simple transposition shuffles of minimum length.

For $1\leq a<n$ denote the \df{simple} transposition $t(a)=t(a,a+1)\in S_n$.  We call a sequence $(a_i)_{i=1}^\ell$ a \df{reduced word} of order $n$ if $\ell={n\choose 2}$ and if the (deterministic) composition $t(a_1)\cdots t(a_\ell)$ equals the reverse permutation $\rho:=[n,\ldots,1]$.  (It is easily verified that $n \choose 2$ is the minimum number of simple transpositions whose composition is $\rho$, and that $\rho$ is the unique permutation in $S_n$ for which this minimum is largest.  Reduced words are extensively studied -- see the later background discussion.)
We construct a simple transposition shuffle from each reduced word as follows.

\begin{construction}[Simple transposition shuffles]\label{adj}
Let $(a_i)_{i=1}^\ell$ be any reduced word.
Write $\sigma_j:=t(a_1)\cdots t(a_j)$ for the composition of the first $j$ transpositions, and let
$$(u_j,v_j):=\bigl(\sigma_{j-1}(a_j), \sigma_{j-1}(a_j+1)\bigr)
$$
be the elements transposed at step $j$.  Let
$\SS=(a_i,a_i+1,p_i)_{i=1}^\ell$ where
\begin{equation*}p_i=\frac{v_i-u_i}{v_i-u_i+1}.%\label{fraction}
\tag*{\dia}
\end{equation*}
\end{construction}

\begin{thm}\label{main}
The minimum length of a simple transposition shuffle of order $n$ is $\ell={n\choose 2}$. If $(a_i)_{i=1}^\ell$ is any reduced word of order $n$, then $\SS$ as defined in \cref{adj} above is a simple transposition shuffle.  Moreover, every simple transposition shuffle of minimum length arises in this way.
\end{thm}

\begin{figure}
\centering
{}\hfill
\begin{tikzpicture}[thick,scale=0.75]
\def\n{5}; \def\l{10};
\def\s{.7};
\foreach \x in {1,...,\n} \draw (\x,-.5*\s)--(\x,-\l*\s-.5*\s);
\foreach \x in {1,...,\n} \node at (\x,0) {$\x$};
\foreach[count=\i] \a/\b/\p in %
{1/2/{1/2}, 2/3/{2/3}, 3/4/{3/4}, 4/5/{4/5}, %
1/2/{1/2}, 2/3/{2/3}, 3/4/{3/4}, %
1/2/{1/2}, 2/3/{2/3}, 1/2/{1/2}} %
{
\node[fill=black,shape=circle,inner sep=0.75mm] (a) at (\a,-\s*\i) {};
\node[fill=black,shape=circle,inner sep=0.75mm] (b) at (\b,-\s*\i) {};
\draw (a)--(b);
\node at (\n+1,-\s*\i) {$\p$};
};
\end{tikzpicture}
\hfill
\begin{tikzpicture}[thick,scale=0.75]
\def\n{5}; \def\l{10};
\def\s{.7};
\foreach \x in {1,...,\n} \draw (\x,-.5*\s)--(\x,-\l*\s-.5*\s);
\foreach \x in {1,...,\n} \node at (\x,0) {$\x$};
\foreach[count=\i] \a/\b/\p in %
{1/2/{1/2}, 3/4/{1/2}, 2/3/{3/4}, 4/5/{2/3}, %
1/2/{2/3}, 3/4/{4/5}, 2/3/{3/4}, %
4/5/{2/3}, 1/2/{1/2}, 3/4/{1/2}} %
{
\node[fill=black,shape=circle,inner sep=0.75mm] (a) at (\a,-\s*\i) {};
\node[fill=black,shape=circle,inner sep=0.75mm] (b) at (\b,-\s*\i) {};
\draw (a)--(b);
\node at (\n+1,-\s*\i) {$\p$};
};
\end{tikzpicture}
\hfill {}

 \medskip

{} \hfill
\begin{tikzpicture}[thick,scale=0.75]
\def\n{5}; \def\l{10};
\def\s{.7};
\foreach \x in {1,...,\n} \draw (\x,-.5*\s)--(\x,-\l*\s-.5*\s);
\foreach \x in {1,...,\n} \node at (\x,0) {$\x$};
\foreach[count=\i] \a/\b/\p in %
{1/2/{1/2}, 2/5/{2/3}, 3/5/{1/4}, 4/5/{1/5}, %
2/3/{1/2}, 1/4/{1/2}, 3/4/{1/2}, %
1/2/{1/2}, 2/3/{2/3}, 1/2/{1/2}} %
{
\node[fill=black,shape=circle,inner sep=0.75mm] (a) at (\a,-\s*\i) {};
\node[fill=black,shape=circle,inner sep=0.75mm] (b) at (\b,-\s*\i) {};
\draw (a)--(b);
\node at (\n+1,-\s*\i) {$\p$};
};
\end{tikzpicture}
\hfill
\begin{tikzpicture}[thick,scale=0.75]
\def\n{5}; \def\l{10};
\def\s{.7};
\foreach \x in {1,...,\n} \draw (\x,-.5*\s)--(\x,-\l*\s-.5*\s);
\foreach \x in {1,...,\n} \node at (\x,0) {$\x$};
\foreach[count=\i] \a/\b/\p in %
{1/2/{1/2}, 3/4/{1/2}, 4/5/{2/3}, 3/4/{1/2}, %
1/3/{(6-\surd 6)/10}, 2/4/{(6+\surd 6)/10}, %
1/2/{1/2}, 3/4/{1/2}, 4/5/{2/3}, 3/4/{1/2}} %
{
\node[fill=black,shape=circle,inner sep=0.75mm] (a) at (\a,-\s*\i) {};
\node[fill=black,shape=circle,inner sep=0.75mm] (b) at (\b,-\s*\i) {};
\draw (a)--(b);
\node at (\n+1.7,-\s*\i) {$\p$};
};
\end{tikzpicture}
\hfill {}
\caption{Examples of transposition shuffles or order $5$ based on \cref{adj} (top left and top right), \cref{sweep} (bottom left), and \cref{divide} (bottom right). A lazy transposition is shown as a horizontal line connecting two positions $a_i$ and $b_i$, with the probability $p_i$ of transposing them given to the right.} \label{examples}
\end{figure}
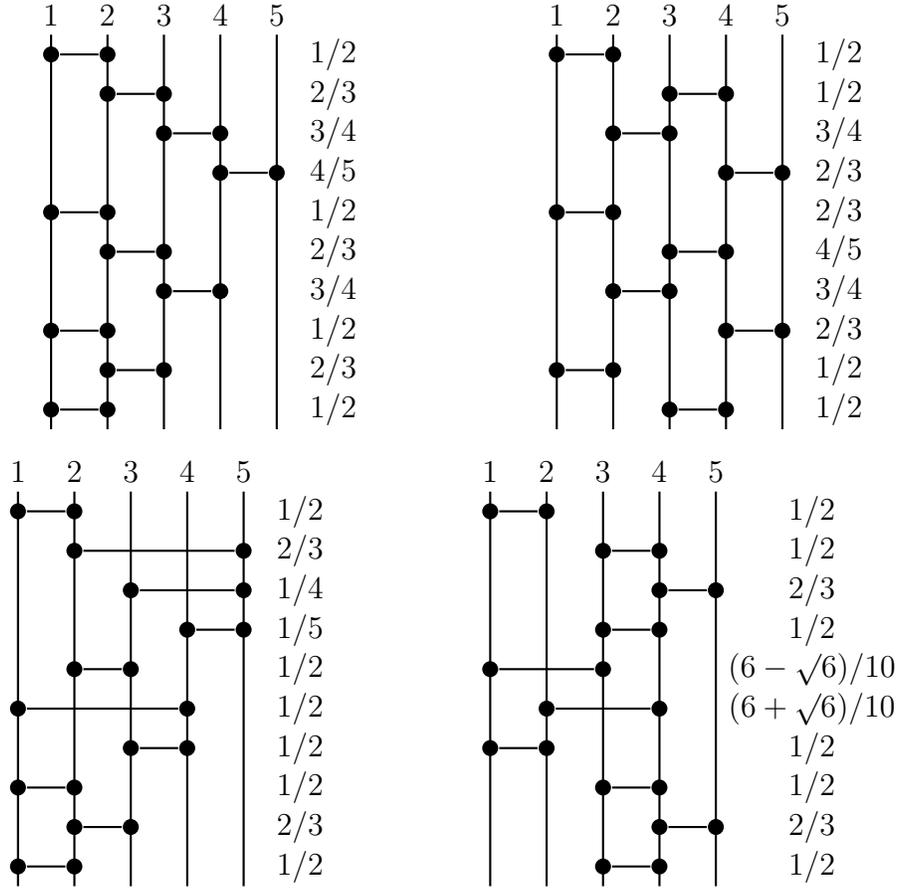

One consequence of \cref{main} is that in any minimum-length simple transposition shuffle, the multiset of probabilities $(p_i)_{i=1}^{\ell}$ comprises exactly $n-1$ $\tfrac12$'s, $n-2$  $\tfrac23$'s, \dots, and one $(n-1)/n$.  See \cref{examples} for examples.

Turning to general (non-simple) transposition shuffles, we will describe several other constructions below, all of length exactly $n \choose 2$, including some with rational probabilities $p_i$ that are not of the form $d/(d+1)$ for integer $d$, and others with irrational probabilities.

For $n\geq 3$ it is not possible for all the probabilities $p_i$ to equal $\tfrac12$, since then the probability of each permutation would be a dyadic rational rather than $1/n!$.  However, we will show that $\tfrac12$ must appear rather frequently.

\begin{thm}\label{half}
In any transposition shuffle $\SS=(a_i,b_i,p_i)_{i=1}^\ell$ of order $n$, we have $$\#\{i: p_i=\tfrac12\}\geq n-1.$$ If the length $\ell$ equals the (in general unknown) minimum $L_n$ then $p_1=p_{\ell}=\tfrac12$.
\end{thm}

\cref{main} implies that in a \emph{simple} transposition shuffle of minimum length, the sequence of probabilities $(p_i)_{i=1}^\ell$ cannot be altered to give another transposition shuffle.  In the general case the following weaker statement holds.

\begin{prop}\label{rigid}
In a transposition shuffle $\SS=(a_i,b_i,p_i)_{i=1}^\ell$ of order $n$ and length $L_n$, the probabilities are rigid in the sense that no single $p_i$ may be altered to give another transposition shuffle.
\end{prop}

\subsection*{Additional constructions}
Next we describe the promised further constructions, together with brief explanations of their correctness and properties.  Also see \cref{examples}.

\begin{construction}[Sweeping]\label{sweep}
We first note the following obvious inductive scheme for constructing transposition shuffles.  Fix $n$.  Call a sequence of parameters $(a_i,b_i,p_i)_{i=1}^k$ a \df{sweep} (of order $n$ and length $k$) if the composition of independent lazy transpositions of $S_n$ with these parameters, $\pi$ say, has the property that its last element $\pi(n)$ is uniformly distributed on $1,\ldots,n$.  The concatenation of any sweep of order $n$ followed by any transposition shuffle of order $n-1$ clearly gives a transposition shuffle of order $n$: the sweep randomizes the final element, then the shuffle shuffles the other elements.

One sweep of length $n-1$ is clearly
$$\bigl(1,2,\tfrac12\bigr),\bigl(2,3,\tfrac23\bigr),\ldots,\bigl(n-1,n,\tfrac{n-1}{n}\bigr).$$
Applying the inductive construction of the previous paragraph to this example gives
a simple transposition shuffle which is a special case of \cref{adj}.
Another straightforward sweep is
$$\bigl(1,n,\tfrac12\bigr),\bigl(2,n,\tfrac13\bigr),\ldots,\bigl(n-1,n,\tfrac1n\bigr).$$

Here is an inductive construction of sweeps, also of length $n-1$, generalizing the last example.  Fix a partition of $\{1,\ldots,n-1\}$ into non-empty sets $D_1,\ldots,D_r$, denote their sizes $m_j:=\#D_j$, and fix an element $d_j\in D_j$ of each.  Apply any sweep of order $m_1$ to $D_1$ in such a way that element $\pi(d_1)$ of the resulting permutation $\pi$ is the one that is uniform on $D_1$.  (Formally, modify the sweep by mapping the parameters $(a_i,b_i,p_i)$ to $(\delta(a_i),\delta(b_i),p_i)$ for each $i$, where $\delta$ is a bijection from $\{1,\ldots,m_1\}$ to $D_1$ mapping $m_1$ to $d_1$.)  Then do similarly for each of $D_2,\ldots,D_r$.  Finally apply a sequence of lazy transpositions with parameters
$$\bigl(d_1,n,\tfrac{m_1}{1+m_1}\bigr),\; \bigl(d_2,n,\tfrac{m_2}{1+m_1+m_2}\bigr),\;
\bigl(d_3,n,\tfrac{m_3}{1+m_1+m_2+m_3}\bigr),\; \ldots,\; \bigl(d_r,n,\tfrac{m_r}{n}\bigr).$$
This ensures that the probability that some element of $D_j$ ends up in location $n$ is $m_j/n$, as required. \dia
\end{construction}

Since the sweeps constructed above all have length $n-1$, the resulting shuffles have length $\sum_{i=1}^{n-1} i={n \choose 2}$.  The length of a sweep of order $n$ must be at least $n-1$, since the graph on vertices $1,\ldots,n$ with edges $\{(a_i,b_i):i=1,\ldots,k\}$ needs to be connected.  So the construction cannot help us to get transposition shuffles of length less than ${n\choose 2}$.

The probabilities $p_i$ that result from \cref{sweep} are all rational, but (unlike those of \cref{adj}) need not be of the form $d/(d+1)$ for integer $d$.  The construction can also give shuffles with $\#\{i: p_i=\tfrac12\}$ strictly greater than $n-1$.

\begin{construction}[Divide and conquer]
\label{divide}
Fix $n$ and let $h=\lfloor n/2 \rfloor$.  Call the integers $1,\ldots,h$ \df{light} and $h+1,\ldots,n$ \df{heavy}.  First apply any transposition shuffle of order $h$, to shuffle the light elements.  Follow this with any transposition shuffle of order $n-h=\lceil n/2\rceil$ on the heavy positions, to shuffle the heavy elements. (Formally, modify each lazy transposition by replacing parameters $(a_i,b_i,p_i)$ with $(a_i+h,b_i+h,p_i)$, and append these triples to the previous list.)

Now apply a sequence of lazy transpositions with parameters
$$(j, j+h, q_j), \quad j=1,\ldots,h,$$
where the probabilities $q_j$ are chosen so that the sum of $j$ independent Bernoulli random variables with parameters $q_1,\ldots,q_h$ is equal in law
to $\#\{i>h:\pi(i)\leq h\}$ where $\pi$ is a uniformly random permutation of $S_n$ -- this is a hypergeometric distribution.  The fact that this is possible is proved in \cite{hypergeom}. (Indeed, the analogous fact holds for a general hypergeometric distribution.  This amounts to the fact that hypergeometric distributions are strongly Rayleigh -- see \cite{strongly-rayleigh}).

At this point, the light and heavy elements are both shuffled, and the number of light elements in heavy positions has the correct distribution.  To complete the construction, we again apply any transposition shuffle to the light positions and apply any transposition shuffle to the heavy positions (as at the start).  This ensures that the locations of the light and heavy elements are shuffled.
\dia
\end{construction}

If the order-$h$ and order-$(n-h)$ shuffles used in \cref{divide} have lengths $\ell_h$ and $\ell_{n-h}$ respectively then the resulting transposition shuffle has length
$$2 \ell_h +2 \ell_{n-h} + h.$$
In particular, if $\ell_h={h \choose 2}$ and $\ell_{n-h}={n-h \choose 2}$ then this is exactly ${n \choose 2}$, so again the construction is of no help in beating this threshold.  On the other hand if it were known that $L_n<{n \choose 2}$ for some fixed $n$ then using \cref{sweep,divide} we could deduce that $L_n\leq (1-\epsilon){n \choose 2}$ for some $\epsilon>0$ and all sufficiently large $n$.

The probabilities $q_j$ in \cref{divide} are in general irrational (but algebraic).  The construction also gives examples of transposition shuffles of length ${n \choose 2}$ in which \emph{two} of the probabilities $p_i$ may be simultaneously altered to give another transposition shuffle (compare \cref{main,rigid}).  For example, two distinct $q_j$ can be exchanged.

\subsection*{Background}
Reduced words have been studied in depth.  For example, it is known \cite{stanley} that the number of reduced words of order $n$ is
$$
\frac{{n \choose 2}!}{1^{n-1} \, 3^{n-2} \, 5^{n-3} \cdots (2n-1)^1}
$$
and that they are in bijection with Young tableau in a certain class \cite{edelman-greene,hamaker-young}.  The uniformly random reduced word of order $n$ has remarkable structure and properties \cite{sort,dauvergne,dauvergne-virag}.

The term reduced word typically refers to a minimum-length sequence of simple transpositions whose composition is an arbitrary specified permutation \cite{saga} (not just $\rho$), and the concept can be extended to general Coxeter (and other) groups.  In \cite{sort} and elsewhere reduced words are referred to as sorting networks (see below).

\sloppypar
Closely related to transposition shuffles are permutation networks and sorting networks.
A sequence $(a_i,b_i)_{i=1}^\ell$ is a \df{permutation network} of order $n$ if for every permutation $\pi\in S_n$ there is some subsequence $j(1),\ldots,j(r)$ of $1,\ldots,\ell$ such that $t(a_{j(1)},b_{j(1)})\cdots t(a_{j(r)},b_{j(r)}) =\pi$.  %(We can interpret a permutation network as a sequence of ``switches'', each of %which can be set to either transpose two elements or do nothing, in such a way %that any permutation is attainable.)
Clearly if $(a_i,b_i,p_i)_{i=1}^\ell$ is a transposition shuffle then $(a_i,b_i)_{i=1}^\ell$ must be a permutation network.
Define the \df{sort} operator $s(a,b)$ by ${x\cdot s(a,b):=x'}$ where for a sequence $x=(x_1,\ldots,x_n)$ the sequence $x'$ agrees with $x$ except that $x'_a=\min(x_a,x_b)$ and $x'_b=\max(x_a,x_b)$.  A \df{sorting network} is a sequence $(a_i,b_i)_{i=1}^\ell$ such that for every permutation $\pi\in S_n$ we have $\pi \cdot s(a_1,b_1)\cdots s(a_\ell,b_\ell)=\id$.  Every sorting network is a permutation network.

There are permutation networks of order $n$ and length asymptotic to $n \log_2 n$ as $n\to \infty$ \cite{waksman}, asymptotically matching the obvious lower bound $\lceil \log_2 n!\rceil \sim n \log_2 n$.  There are sorting networks of length $O(n\log n)$ \cite{aks}, but known constructions are quite indirect and complex, with impractically large constants in the $O$ notation; on the other hand there are straightforward constructions of length $O(n \log^2 n)$ with reasonable constants \cite{batcher}.  Can these networks be turned into transposition shuffles?

Restricting attention to simple transpositions, $(a_i,a_i+1)_{i=1}^\ell$ is a sorting network if and only if it is a permutation network, and moreover this is equivalent to the condition $t(a_1)\cdots t(a_\ell)=\rho$;  see e.g.\ \cite[5.3.4]{knuthvol3}.  Thus, the minimum length of a simple permutation network (or sorting network) is ${n\choose 2}$, and the minimal examples coincide with reduced words as defined earlier.

Sorting networks have applications in distributed or hardware-optimized systems such as graphics processor units.  Transposition shuffles also appear natural for applications, since the ability to permute objects uniformly is useful for privacy or security as well as for games of chance.

To our knowledge transposition shuffles have not been considered before. The problem of ``square roots'' of uniform measure addressed in \cite{square-roots} is somewhat related, while \emph{asymptotic} shuffling under various random transposition models has been studied extensively -- see e.g.\ \cite{lpw} for a comprehensive treatment and \cite{extra} for a specific model close to the one considered here.

\section{Simple transpositions}

We divide the proof of \cref{main} into two parts.  First we show that \cref{adj} works; then we show that it exhausts the possibilities.  We will use several standard properties of reduced words.  For a reduced word $(a_i)_{i=1}^\ell$ recall that we write $\sigma_{j}=t(a_1)\cdots t(a_{j})$ for the (deterministic) permutation after $j$ steps, so that in particular $\sigma_0=\id$ and $\sigma_\ell=\rho$.  On the other hand we write $\pi_j=T_1\cdots T_j$ for the (random) composition of the first $j$ lazy transpositions, so that in a transposition shuffle $\pi_\ell$ is uniform on $S_n$.

A reduced word $(a_i)_{i=1}^\ell$ may be transformed into another via \df{moves} of the following types.
\begin{ilist}
\item \df{Commuting move:} if two consecutive elements $a_{j},a_{j+1}$ satisfy $|a_j-a_{j+1}|\geq 2$, exchange them to get the word $(a_1,\ldots,a_{j+1},a_j,\ldots,a_\ell)$.
\item \df{Braid move:} if three consecutive elements $(a_{j},a_{j+1},a_{j+2})$ are of the form $(k,k+1,k)$, replace them with $(k+1,k,k+1)$, or vice versa.
\end{ilist}

\begin{prop}[Tits, \cite{tits}]\label{tits}
Any reduced word may be transformed into any other via a sequence of moves of types (i) and (ii).
\end{prop}

\cref{tits} is a special case of a more general result \cite{tits}, which applies to reduced words of an arbitrary permutation (not just $\rho$), and to general Coxeter groups.  We remark that the result would be essentially obvious if we in addition allowed moves of the form $(k,k)\leftrightarrow()$, which change the length of the word.

\begin{figure}
\centering
{}\hfill
\begin{tikzpicture}[thick,scale=0.75]
\def\n{3}; \def\l{3};
\def\s{.9};
\foreach \x in {1,...,\n} \draw (\x,-.5*\s)--(\x,-\l*\s-.5*\s);
\foreach[count=\i] \a/\b/\p in %
{1/2/{1/2}, 2/3/{4/5}, 1/2/{3/4}} %
{
\node[fill=black,shape=circle,inner sep=0.75mm] (a) at (\a,-\s*\i) {};
\node[fill=black,shape=circle,inner sep=0.75mm] (b) at (\b,-\s*\i) {};
\draw (a)--(b);
\node at (\n+1,-\s*\i) {$\p$};
};
\end{tikzpicture}
\hfill
\begin{tikzpicture}[thick,scale=0.75]
\def\n{3}; \def\l{3};
\def\s{.9};
\foreach \x in {1,...,\n} \draw (\x,-.5*\s)--(\x,-\l*\s-.5*\s);
\foreach[count=\i] \a/\b/\p in %
{2/3/{3/4}, 1/2/{4/5}, 2/3/{1/2}} %
{
\node[fill=black,shape=circle,inner sep=0.75mm] (a) at (\a,-\s*\i) {};
\node[fill=black,shape=circle,inner sep=0.75mm] (b) at (\b,-\s*\i) {};
\draw (a)--(b);
\node at (\n+1,-\s*\i) {$\p$};
};
\end{tikzpicture}
\hfill {}
\caption{A braid move.  The two sequences of independent lazy transpositions have the same effect.}
\label{braid-ex}
\end{figure}
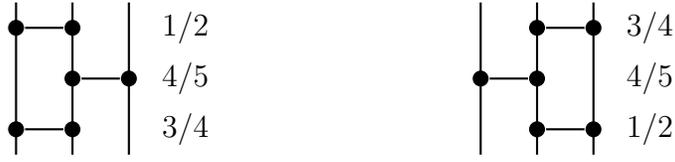

We next address how to transform probabilities under braid moves.  See \cref{braid-ex} for an example.
\begin{samepage}
\begin{lemma}\label{braid}
Let $n=3$, let $(T_i)_{i=1}^3$ be independent lazy transpositions with respective parameters $(1,2,p),(2,3,q),(1,2,r)$, and let $(T'_i)_{i=1}^3$ be independent lazy transpositions with respective parameters $(2,3,p'),(1,2,q'),(2,3,r')$.  Given $p,q,r\in(0,1)$, it is possible to choose $p',q',r'\in(0,1)$ so that the compositions $T_1T_2T_3$ and $T'_1T'_2T'_3$ are equal in law if and only if $$\frac{p}{1-p}+\frac{r}{1-r}=\frac{q}{1-q},$$ in which case the unique choice is $(p',q',r')=(r,q,p)$.
\end{lemma}
\end{samepage}

We remark that the requirement that the probabilities lie in $(0,1)$ is needed -- otherwise there are further (degenerate) solutions such as $(p,q,r)=(1,0,0)$, $(p',q',r')=(0,1,0)$.

\begin{proof}[Proof of \cref{braid}]
It is convenient to use the parameters $P=p/(1-p)$, $Q=q/(1-q)$ etc.  The two compositions are equal in law if and only if they assign equal probabilities to each of the $3!$ permutations.  This yields the equations
\begin{align*}
PQR&=P'Q'R'&
PQ&=Q'R'&
QR&=P'Q'\\
1+PR&=1+P'R'&
Q&=P'+R'&
P+R&=Q',
\end{align*}
which are equivalent to the conditions $P+R=Q$ and $(P',Q',R')=(R,Q,P)$.
\end{proof}

\begin{proof}[Proof of \cref{main} -- correctness of construction]
Here we prove that for any reduced word, \cref{adj} gives a transposition shuffle.
We first check this in an easy case, the ``bubble sort'' word.  Let $\ell={n\choose 2}$ and take:
\begin{alignat*}{5}
(a_i)_{i=1}^{\ell} &= \bigl(1,2,\ldots,n-1,\ &&1,2,\ldots,n-2,\ &&\ldots\  &&2,1,\ &&1\bigr).\\
\intertext{The corresponding probabilities according to \cref{adj} are:}
(p_i)_{i=1}^{\ell} &= \bigl(\tfrac12,\tfrac23,\ldots,\tfrac{n-1}{n},\ &&\tfrac12,\tfrac23,\ldots,\tfrac{n-2}{n-1},\ &&\ldots\ &&\tfrac12,\tfrac23,\ &&\tfrac12\bigr).
\end{alignat*}
We check that $\SS=(a_i,a_i+1,p_i)_{i=1}^\ell$ is a transposition shuffle.  Indeed,
this is a special case of \cref{sweep}: the first $n-1$ steps form a sweep, so the permutation $\pi_{n-1}$ has uniformly random last element $\pi_{n-1}(n)$.  The remaining sequence of parameters agrees with the entire sequence for $n-1$, so by induction, their composition is a uniformly random permutation of elements $1,\ldots,n-1$, concluding the argument.

Now we apply \cref{tits}.  Suppose that $(a_i)_{i=1}^{\ell}$ and $(a'_i)_{i=1}^{\ell}$ are reduced words that are related by a single move, and let $(p_i)_{i=1}^{\ell}$ and $(p'_i)_{i=1}^{\ell}$ be the corresponding probabilities given by \cref{adj} in each case.  It suffices to show that if $(a_i,a_i+1,p_i)_{i=1}^\ell$ is a transposition shuffle then so is $(a'_i,a'_i+1,p'_i)_{i=1}^\ell$.  This clearly holds in the case of a commuting move.

For a braid move we will use \cref{braid}.  Suppose without loss of generality that the move replaces $(a_j,a_{j+1},a_{j+2})=(k,k+1,k)$ with $(a'_j,a'_{j+1},a'_{j+2})=(k+1,k,k+1)$.  Write $(u,v,w)=(\sigma_{j-1}(k),\sigma_{j-1}(k+1),\sigma_{j-1}(k+2))$ for the three (deterministic) elements involved, which satisfy $u<v<w$.  Let $(p,q,r)=(p_j,p_{j+1},p_{j+2})$ and $(p',q',r')=(p'_j,p'_{j+1},p'_{j+2})$ be the three probabilities before and after the move, and write $P=p/(1-p)$ so that $p=P/(1+P)$, etc. Then the formula for the probabilities in \cref{adj} gives
\begin{align*}
(P,Q,R)&=\bigl(v-u,w-u,w-v\bigr)\\
(P',Q',R')&=\bigl(w-v,w-u,v-u\bigr).
\end{align*}
These values satisfy the conditions of \cref{braid}.
\end{proof}

We now prepare for the uniqueness part of the proof of \cref{main}.  We parameterize the probability space as follows.  Given $(a_i,b_i,p_i)_{i=1}^\ell$ as usual, let $\omega=(\omega_i)_{i=1}^\ell$ be independent $\{0,1\}$-valued random variables with $\P(\omega_i=1)=p_i$, and let the lazy transposition $T_i$ equal $t(a_i,b_i)$ if and only if $\omega_i=1$ (and otherwise equal $\id$).  Thus, for given $(a_i,b_i)_{i=1}^\ell$, the permutations $\pi_0,\ldots,\pi_\ell$ are deterministic functions of $\omega$.  They take values $\sigma_0,\ldots,\sigma_\ell$ when $\omega$ is the all-$1$ vector.

\begin{figure}
{} \hfill
\begin{tikzpicture}[thick,scale=0.75]
\def\n{5}; \def\l{10};
\def\s{.7};
\foreach \x in {1,...,\n} \draw (\x,-.5*\s)--(\x,-\l*\s-.5*\s);
\node[text=red] at (1,0) {$1$};
\node[text=red] at (5,-7.7) {$1$};
\foreach[count=\i] \a/\b/\p in %
{1/2/{1/2}, 3/4/{1/2}, 2/3/{3/4}, 4/5/{2/3}, %
1/2/{2/3}, 3/4/{4/5}, 2/3/{3/4}, %
4/5/{2/3}, 1/2/{1/2}, 3/4/{1/2}} %
{
\node[fill=black,shape=circle,inner sep=0.75mm] (a) at (\a,-\s*\i) {};
\node[fill=black,shape=circle,inner sep=0.75mm] (b) at (\b,-\s*\i) {};
\draw (a)--(b);
};
\draw[line width=3,red] (1,-.35)--(1,-.7)--(2,-.7)--(2,-2.1)--%
(3,-2.1)--(3,-4.2)--(4,-4.2)--(4,-5.6)--(5,-5.6)--(5,-7.35);
\end{tikzpicture}
\hfill
\begin{tikzpicture}[thick,scale=0.75]
\def\n{4}; \def\l{10};
\def\s{.7};
\foreach \x in {1,...,\n} \draw (\x,-.5*\s)--(\x,-\l*\s-.5*\s);
\node[text=red] at (1,0) {\phantom{$1$}};
\node[text=red] at (5,-7.7) {\phantom{$1$}};
\foreach \a/\b/\i in %
{ 2/3/{2}, 3/4/{4}, %
1/2/{5},  2/3/{7}, %
 1/2/{9}, 3/4/{10}} %
{
\node[fill=black,shape=circle,inner sep=0.75mm] (a) at (\a,-\s*\i) {};
\node[fill=black,shape=circle,inner sep=0.75mm] (b) at (\b,-\s*\i) {};
\draw (a)--(b);
};
\end{tikzpicture}
\hfill {}
\caption{An illustration of \cref{traj} and its application in the proof of \cref{main}. Left: the only way element $1$ can reach location~$n$.  Right: the lower-order reduced word obtained by deleting its trajectory.}
\label{delete}
\end{figure}
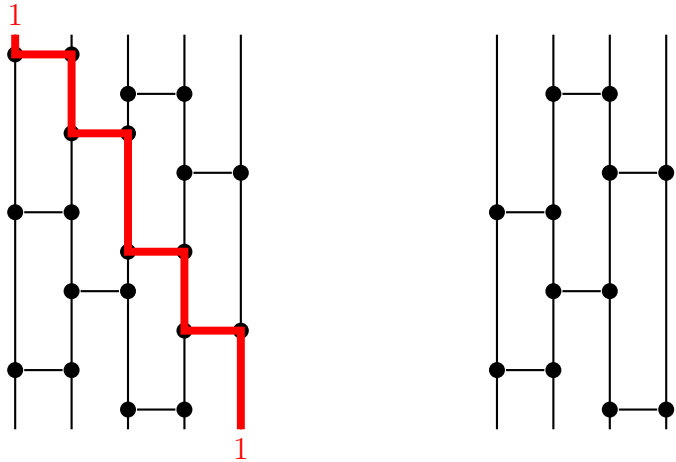

The next lemma characterizes the ways in which element $1$ can reach position $n$ in the final permutation $\pi_\ell$.  See \cref{delete}.

\begin{samepage}
\begin{lemma}\label{traj}
For any reduced word $(a_i)_{i=1}^\ell$ of order $n$ there exists a fixed set $H\subseteq\{1,\ldots,\ell\}$ with $|H|=n-1$ such that $\pi_\ell(n)=1$ if and only if $\omega_h=1$ for all $h\in H$.  Moreover, in that case the trajectory of element $1$ satisfies (and is determined by): %$\pi_i^{-1}(1)=1+|H\cap [1,i]|$
$\pi_i^{-1}(1)-\pi_{i-1}^{-1}(1)=\ind[i\in H]$ for all $0< i \leq \ell$.
\end{lemma}
\end{samepage}

\begin{proof}
First note that in the determinisitic sequence of permutations $\sigma_0,\ldots,\sigma_\ell$, each of the $n \choose 2$ pairs of elements of $\{1,\ldots,n\}$ swaps exactly once.  In particular, element $1$ swaps with every other element exactly once, and it must move from left to right when it does so.  Let
\begin{align*}
H&=\bigl\{h: \sigma_h^{-1}(1)-\sigma_{h-1}^{-1}(1)=1\bigr\}\\
&=\bigl\{h: \sigma_{h-1}^{-1}(1)=a_h\bigr\}
\end{align*}
be the set of times when these swaps occur.  Note that no other transpositions are incident to the trajectory of $1$; that is,
\begin{equation}
\sigma_{i-1}^{-1}(1)\notin\{a_i,a_i+1\}\quad\text{for }i\notin H.
\label{inc}
\end{equation}

It follows immediately that $\pi_i^{-1}(1)=\sigma_i^{-1}(1)$ for all $i$ if and only if $\omega_h=1$ for all $h\in H$.  It remains to show that there is no other possible trajectory via which element $1$ can end at position $n$.  Suppose on the contrary that $\omega$ is such that $\pi_\ell(1)=n$ but $\pi_i^{-1}(1)\neq\sigma_i^{-1}(1)$ for some $i$.  If $\pi_i^{-1}(1)<\sigma_i^{-1}(1)$ for some $i$, consider the largest $i$ for which this holds.  Then we must have
$$1+\pi_i^{-1}(1)=\sigma_i^{-1}(1)=\pi_{i+1}^{-1}(1)=\sigma_{i+1}^{-1}(1)=k,\text{ say}.$$
But this implies that $a_{i+1}=k-1$ (and $\omega_{i+1}=1$), giving a contradiction to \eqref{inc} and the definition of $H$.  On the other hand if $\pi_i^{-1}(1)>\sigma_i^{-1}(1)$ for some $i$, considering the smallest such $i$ leads similarly to a contradiction.
\end{proof}

\begin{proof}[Proof of \cref{main} -- uniqueness]
It is clear that no simple transposition shuffle can have length less than $n\choose 2$, since it would be incapable of producing the reverse permutation $\rho$ (in which every pair of elements is reversed).

It remains to show uniqueness: for any reduced word $(a_i)_{i=1}^\ell$ there is at most one sequence of probabilities $(p_i)_{i=1}^\ell$ for which $(a_i,a_{i}+1,p_i)_{i=1}^\ell$ is a transposition shuffle.  We prove this statement by induction on the order $n$.  It is clearly true for $n\leq 2$.

Fix the reduced word, and let $H$ be the set from \cref{traj}.  By that lemma, for any transposition shuffle we must have
$$\frac1n=\P\bigl(\pi_\ell(n)=1\bigr)=\P\bigl(\omega_h=1\,\forall h\in H\bigr)=
\prod_{h\in H} p_h,$$
so this product of $p_h$'s is uniquely determined.  Moreover, conditional on the event $\pi_\ell(n)=1$, the remaining elements $[\pi_\ell(1),\ldots,\pi_\ell(n-1)]$ should form a uniformly random permutation of $2,\ldots,n$.

To make use of this last fact we delete the trajectory of element $1$ from the reduced word to get a lower-order word.  More precisely, define $(c_i)_{i=1}^\ell$ by
$$c_i:=\begin{cases}
a_i, &a_i<\sigma_{i-1}^{-1}(1);\\
a_i-1, &a_i>\sigma_{i-1}^{-1}(1)+1;\\
\infty, &i\in H,
\end{cases}$$
where we use $\infty$ as a dummy symbol.
Then the subsequence $(a'_i)_{i=1}^{\ell'}:=(c_i:i\notin H)$ obtained by deleting all occurrences of $\infty$ is a reduced word of order $n-1$ (and length $\ell':={n\choose 2}-(n-1)={n-1 \choose 2}$).  See \cref{delete}.
Now consider any $\omega\in\{0,1\}^\ell$ that satisfies $\omega_h=1$ for all $h\in H$, and define the subsequence $\omega':=(\omega_i:i\notin H)\in\{0,1\}^{\ell'}$.  Let $\pi'_{\ell'}$ be the final permutation of an order-$(n-1)$ simple transposition shuffle with word $(a'_i)_{i=1}^{\ell'}$ at the element $\omega'$ of its probability space.  Then the final permutation under the original shuffle at $\omega$ is
$$\pi_\ell=\bigl[\pi'_{\ell'}(1)+1,\ldots,\pi'_{\ell'}(n-1)+1,1\bigr].$$
Therefore, by the induction hypothesis, there is at most one possible choice of the vector of probabilities $(p_i)_{i\notin H}$ that results in the correct conditional law of the permutation $\pi_\ell$ given $\pi_\ell(n)=1$.

Now, by symmetry, we can apply the same argument to the set
$$\widehat{H}=\bigl\{h: \sigma_h^{-1}(n)-\sigma_{h-1}^{-1}(n)=-1\bigr\}$$
of times when element $n$ moves, to deduce that there is also at most one choice for the vector of probabilities $(p_i)_{i\notin \widehat{H}}$.  Now, $H\cap\widehat H$ has exactly one element: it is the unique time $k$ at which elements $1$ and $n$ swap in $\sigma_0,\ldots,\sigma_{\ell}$.  Hence there is at most one choice for $(p_i)_{i\neq k}$.  But since $\prod_{h\in H} p_h$ is determined, there is at most one choice for $p_k$ also.
\end{proof}

\section{General Transpositions}

\begin{proof}[Proof of \cref{rigid}]
We fix the parameters $(a_i,b_i)_{i=1}^\ell$ and consider dependence of the law of the final permutation $\pi_\ell$ on the probabilites $(p_i)_{i=1}^\ell$.  For any given permutation $\alpha\in S_n$ we have
$$\P(\pi_\ell=\alpha)=\sum_{\omega\in S_\alpha}\prod_{i:\omega_i=1} p_i \prod_{i:\omega_i=0} (1-p_i)$$
for some set $S_\alpha\subseteq\{0,1\}^\ell$.
Suppose we vary one probability $p_j$ while fixing the others.  Then the dependence is affine:
$$\P(\pi_\ell=\alpha)=A+B p_j,$$
where the constants $A$ and $B$ depend on $j$, $\alpha$ and $(p_i:i\neq j)$.
Suppose that the choice of probabilities $(p_i)_{i=1}^\ell$ gives a transposition shuffle, and so also do the probabilities obtained by altering $p_j$ (only) to a different value $p'_j\neq p_j$.  Then in particular
$$A+B p_j = A+B p'_j=1/n!,$$
so $B=0$, hence \emph{any} choice of $p''_j\in[0,1]$ will also give $\P(\pi_\ell=\alpha)=1/n!$.  The same argument applies for every permutation $\alpha$, so any choice of $p''_j$ gives a transposition shuffle.  But in particular we can take $p''_j=0$ and remove the $j$th lazy transposition altogether, so $\ell$ was not minimal.
\end{proof}

\begin{proof}[Proof of \cref{half}]
We first prove the statement about the first and last probabilities.  Suppose $\mathcal{S}=(a_i,b_i,p_i)_{i=1}^\ell$ is a minimum-length transposition shuffle.  Appending a lazy transposition with parameters
$(a_\ell,b_\ell,\tfrac12)$
to the sequence clearly gives another transposition shuffle.  But now the last two lazy transpositions can be replaced a single one of parameters
$(a_\ell,b_\ell, \frac12)$.  This contradicts rigidity, \cref{rigid}, unless $p_\ell=\tfrac12$.   Symmetry gives $p_1=\tfrac12$ also.

Now we turn to the claim about the number of occurrences of $\tfrac12$.  To any random permutation $\pi$ of $S_n$ we can associate the $n\times n$ matrix $M(\pi)$ with entries
$$M(\pi)_{i,j}=\P(\pi(i)=j).$$
(In other words, $M(\pi)$ is the expectation of the permutation matrix.)  If $\pi$ and $\tau$ are independent random permutations then $M(\pi\tau)=M(\pi)M(\tau)$.

In a transposition shuffle, $M(\pi_\ell)$ is the matrix with all entries $1/n$, which has rank $1$.  On the other hand, the matrix $M(T)$ of the lazy transposition with parameters $(a,b,p)$ agrees with the identity except in the intersection of rows $a$ and $b$ with columns $a$ and $b$, where it has the form
$$\left(
    \begin{array}{cc}
      1-p & p \\
      p & 1-p \\
    \end{array}
  \right).
$$
Thus $M(T)$ has rank $n$ if $p\neq\tfrac12$ and rank $n-1$ if $p=\tfrac12$.  Sylvester's rank inequality states that $n-\rank(AB)\leq n-\rank(A)+n-\rank(B)$ for $n\times n$ matrices $A,B$, so we deduce that $\{i:p_i=\tfrac12\}\geq n-1$ as required.
\end{proof}

\section*{Acknowledgements}
We thank Swee Hong Chan for helpful comments on an earlier draft, and Viktor Kiss for sharing the results of his computational experiments.

\bibliographystyle{abbrv}
\bibliography{shuf}

\end{document}